\documentclass[12pt,reqno]{amsart}
\usepackage{amssymb, amsmath, hyperref}
\usepackage{pdfsync}
\usepackage{stmaryrd}
\usepackage{bbold}
\usepackage{bbm}
\usepackage{graphics,graphicx,fancyhdr,color}

\newtheorem{df}{Definition}[section]

\newtheorem{thm}[df]{Theorem}
\newtheorem{cor}[df]{Corollary}

\makeatletter \@addtoreset{equation}{section}

\newcommand{\bes}{\begin{displaymath}}
\newcommand{\ees}{\end{displaymath}}
\newcommand{\be}{\begin{equation}}
\newcommand{\ee}{\end{equation}}
\newcommand{\ba}{\begin{eqnarray}}
\newcommand{\ea}{\end{eqnarray}}
\newcommand{\bas}{\begin{eqnarray*}}
\newcommand{\eas}{\end{eqnarray*}}
\newcommand{\@Bbb}[1]{\ensuremath{\mathbb #1}}

\newcommand{\B}{{\@Bbb B}}
\newcommand{\C}{{\@Bbb C}}

\newcommand{\F}{{\@Bbb F}}

\newcommand{\Q}{{\@Bbb Q}}
\newcommand{\bQ}{{\@Bbb Q}}
\newcommand{\N}{{\@Bbb N}}

\newcommand{\bbR}{{\@Bbb R}}
\newcommand{\W}{{\@Bbb W}}

\newcommand{\bbZ}{{\@Bbb Z}}
\newcommand{\bbT}{{\@Bbb T}}

\newcommand{\eps}{\epsilon}

\newcommand{\@s}[1]{\ensuremath{\mathcal #1}}
\newcommand{\cA}{\@s A}
\newcommand{\cB}{\@s B}
\newcommand{\cC}{\@s C}
\newcommand{\cD}{\@s D}
\newcommand{\cE}{\@s E}
\newcommand{\cF}{\@s F}
\newcommand{\cG}{\@s G}
\newcommand{\cH}{\@s H}
\newcommand{\cI}{\@s I}
\newcommand{\cJ}{\@s J}

\newcommand{\cK}{\@s K}
\newcommand{\cL}{\@s L}
\newcommand{\cN}{\@s N}
\newcommand{\cM}{\@s M}
\newcommand{\cO}{\@s O}
\newcommand{\cP}{\@s P}
\newcommand{\cR}{\@s R}
\newcommand{\cS}{\@s S}
\newcommand{\cT}{\@s T}
\newcommand{\cV}{\@s V}
\newcommand{\cW}{\@s W}
\newcommand{\cX}{\@s X}
\newcommand{\cY}{\@s Y}
\newcommand{\cZ}{\@s Z}

\newcommand{\@bm}[1]{\ensuremath{\mathbf #1}}
\newcommand{\bma}{\@bm a}
\newcommand{\bmb}{\@bm b}
\newcommand{\bmc}{\@bm c}
\newcommand{\bmd}{\@bm d}

\newcommand{\bme}{\@bm e}
\newcommand{\bmf}{\@bm f}
\newcommand{\bmg}{\@bm g}
\newcommand{\bmh}{\@bm h}
\newcommand{\bmi}{\@bm i}
\newcommand{\bmj}{\@bm j}
\newcommand{\bmk}{\@bm k}
\newcommand{\bml}{\@bm l}
\newcommand{\bmm}{\@bm m}
\newcommand{\bmn}{\@bm n}
\newcommand{\bmo}{\@bm o}
\newcommand{\bmp}{\@bm p}
\newcommand{\bmq}{\@bm q}
\newcommand{\bmr}{\@bm r}
\newcommand{\bms}{\@bm s}
\newcommand{\bmt}{\@bm t}
\newcommand{\bmu}{\@bm u}
\newcommand{\bmw}{\@bm w}
\newcommand{\bmv}{\@bm v}
\newcommand{\bmx}{\@bm x}
\newcommand{\bx}{\@bm x}

\newcommand{\bmy}{\@bm y}
\newcommand{\bz}{\@bm z}

\newcommand{\by}{\@bm y}

\newcommand{\bmzero}{\@bm 0}

\newcommand{\ga}{\gamma}

\newcommand{\@g}[1]{\ensuremath{\mathfrak #1}}
\newcommand{\gA}{\@g A}
\newcommand{\gD}{\@g D}
\newcommand{\gJ}{\@g J}
\newcommand{\gF}{\@g F}
\newcommand{\gM}{\@g M}
\newcommand{\gR}{\@g R}

\newcommand{\commentout}[1]{{}}

\makeatother

\begin{document}

\title[Infinite hermitian  matrices]{A criterion for  essential self-adjointness of
  a symmetric operator defined by some infinite hermitian matrix with
  unbounded entries}

\author{Tomasz Komorowski}
\begin{abstract}
We shall consider a double infinite, hermitian, complex
entry matrix
$A=[a_{x,y}]_{x,y\in\bbZ}$, with $a_{x,y}^*=a_{y,x}$, $x,y\in\bbZ$.
Assuming that
the matrix is {\em almost of a finite bandwidth}, i.e. 
 there exists an integer $n> 0$ and exponent $\ga\in[0,1)$ such
that 
$
a_{x,x+z}=0$ for all $z>n\langle x\rangle^{\ga}$ and the growth of the
$\ell_1$ norm of a row is slower than $|x|^{1-\ga}$ for $|x|\gg1$,
i.e. $\lim_{|x|\to+\infty}| x|^{\ga-1}\sum_{y}|a_{xy}|=0$ we prove
that the corresponding symmetric operator, defined on compactly supported sequences,  is essentially
self-adjoint in $\ell_2(\bbZ)$. In the case $\ga=0$ (the so called $(nJ)$-matrices) we
prove that there exists $c_*>0$, depending only on $n$, such that  the condition
$\limsup_{|x|\to+\infty}| x|^{-1}\sum_{y}|a_{xy}|\le c_*$ suffices to
conclude essential self-adjointness.
\end{abstract}

\vspace*{-1in}

\maketitle

\section{Introduction}

We shall consider a double infinite, hermitian, complex
entry matrix
$A=[a_{x,y}]_{x,y\in\bbZ}$, with $a_{x,y}^*=a_{y,x}$,
$x,y\in\bbZ$. Here $a^*$ denotes the complex conjugate of 
$a\in\mathbb C$.
We assume furthermore that
the matrix is {\em almost of a finite bandwidth}, i.e. 
 there exists an integer $n\ge 1$ and exponent $\ga\in[0,1)$ such
that 
\begin{equation}
\label{finite-bw}
a_{x,x+z}=0\mbox{ for all }z>n\langle x\rangle^{\ga}.
\end{equation}
Here, for  given $a$ we let
$\langle a\rangle:=(1+|a|^2)^{1/2}.$
With
the help of matrix $A$ we can define a symmetric operator on the subset $c_0(\bbZ)$ of the
complex Hilbert space $\ell_2(\bbZ)$ --  the space consisting of all double infinite
sequences $f=(f_x)$ equipped with the norm
$$
\|f\|_{\ell_2(\bbZ)}:=\left\{\sum_x|f_x|^2\right\}^{1/2}<+\infty.
$$
Here  $c_0(\bbZ)$ is the subspace containing all compactly
supported $f$.  
The operator is given by 
\begin{equation}
\label{052609}
(Af)_x:=\sum_{y}a_{xy}f_y,\quad x\in\bbZ,\quad f\in c_0(\bbZ).
\end{equation}
According to Theorem  4,
p. 102 of \cite{akhiezer},  assumption \eqref{finite-bw} implies that
the operator is closable. Denote its closure by 
$\bar A:D(\bar A)\to\ell_2(\bbZ)$.
 In our
principal result, see Theorem \ref{main} below, we formulate a
sufficient condition, in terms of the growth of  $|a_{xy}|$, see \eqref{042609} below, for the
operator $\bar A$ to be self-adjoint. The above means that  the
deficiency index of $A:c_0(\bbZ)\to\ell_2(\bbZ)$ equals $(0,0)$, see \cite{nagy}.

Note that in the particular case when $\ga=0$ we have
\begin{equation}
\label{finite-bw-1}
a_{x,y}=0\mbox{ for all }|x-y|>n
\end{equation} 
and the definition coincides
with the usual definition of $(nJ)$-matrices, see \cite{smart},
(sometimes also called {\em  finite bandwidth}
matrices). When $n=1$
they are called {\em Jacobi} matrices and play an important role in
the theory of Hamburger moment problem. This case has been well
studied in the literature, see e.g. \cite{kost, krein1, krein2, kogan,
  simon}  and the references contained therein, although also then our results formulated in Corollary \ref{cor-main} and Theorem \ref{cor-main1} below seem to be new.


\bigskip


\bigskip

\section{The statement of the main result}

\bigskip

Since matrix $A=[a_{x,y}]$ is hermitian the operator $\bar A$ is obviously symmetric, i.e.
\begin{equation}
\label{022609}
\langle g,\bar A f\rangle_{\ell_2(\bbZ)}=\langle
\bar Ag,f\rangle_{\ell_2(\bbZ)},\quad f,g\in D(\bar A).
\end{equation}
Here, $\langle \cdot,\cdot\rangle_{\ell_2(\bbZ)}$ denotes  the usual scalar product in $\ell_2(\bbZ)$.
Let
$f\in\ell_2(\bbZ)$ be such that the functional
\begin{equation}
\label{012609}
\varphi(g):=\langle \bar Ag,f\rangle_{\ell_2(\bbZ)},\quad g\in D(\bar A)
\end{equation}
is bounded, i.e. for some $C>0$
\begin{equation}
\label{032609}
|\varphi(g)|\le C\|g\|_{\ell_2},\quad g\in D(\bar A).
\end{equation}
Self-adjointness of $\bar A$ means that any $f$, for which \eqref{032609}
holds, belongs to $D(\bar A)$ and, as a consequence, \eqref{022609}
is in force.

For example, if there exists $M>0$ such that
 $\sum_{y}|a_{xy}|\le M$ for all $x\in\bbZ$  then $\bar A$ is
bounded 
on $\ell_2(\bbZ)$, see Example III.2.3, p. 143 of \cite{kato}, therefore it is self-adjoint. 
Our main result can be stated as follows.
\begin{thm}
\label{main}
Suppose that  for some $\ga\in[0,1)$ the entries of matrix $A$ satisfy
both condition
\eqref{finite-bw} and
\begin{equation}
\label{042609}
\lim_{|x|\to+\infty}\frac{\!\!1}{\langle x\rangle^{1-\ga}}\left(\sum_{y}|a_{xy}|\right)=0.
\end{equation}
Then, operator $A$, given by \eqref{052609}, is 
 essentially selfadjoint on $\ell_2(\bbZ)$.
\end{thm}

Using  the theorem for $\ga=0$ we immediately conclude the following.
\begin{cor}
\label{cor-main}
The conclusion of Theorem \ref{main} holds when  $A$ is a hermitian,
$nJ$-matrix (i.e. \eqref{finite-bw-1} is in force) whose entries satisfy
\begin{equation}
\label{042609a}
\lim_{|x|\to+\infty}\frac{1}{\langle x\rangle}\left(\sum_{y}|a_{xy}|\right)=0.
\end{equation}
\end{cor}

\bigskip

In fact, in the case of
$nJ$-matrices, one can show a little stronger result, relaxing  a bit
assumption \eqref{042609a}.
\begin{thm}
\label{cor-main1}
There exists
$c_*>0$ depending only on $n$ such that the conclusion of Theorem
\ref{main} holds  for any  hermitian $nJ$-matrix $A=[a_{xy}]$ that satisfies 
\begin{equation}
\label{042609b}
\limsup_{|x|\to+\infty}\frac{1}{\langle x\rangle}\left(\sum_{y}|a_{xy}|\right)\le c_*.
\end{equation}
\end{thm}
\noindent{\bf Example.}  
The condition \eqref{042609b} is in some sense optimal. Suppose that
$\delta>1$ is arbitrary. Consider the Jacobi
matrix with entries given by
$$
a_{x,x+z}=\left\{
\begin{array}{ll}
0, & \mbox{if $x\le 0$, or $x+z\le0$, or $z=0$, or $z>2$},\\
&\\
x^{\delta}, & \mbox{if $x>0$ and $z=1$}.
\end{array}
\right.
$$ 
According to Corollary 1, p. 267, of \cite{kost} the index of
deficiency of the respective operator $A:c_0(\bbZ)\to\ell_2(\bbZ)$ equals then $(1,1)$. Therefore
$A$  cannot be essentially self-adjoint, see Corollary 2.2 of \cite{simon}.
\bigskip

\section{Proof of Theorem \ref{main}}

Recall the classical criterion for the essential self-adjointness of a
symmetric operator, see Theorem 3 of Section 33.2
of \cite{lax}, or Corollary 2.2 of \cite{simon}. Adjusted to our settings it reads as follows: suppose that a closed operator $\bar A$ is symmetric and
\begin{equation}
\label{062609}
R(I-i\bar A) =\ell_2(\bbZ)=R(I+i\bar A).
\end{equation}
Then, it  is self-adjoint.

To prove \eqref{062609} we show that for any $g=(g_x)\in c_0(\bbZ)$ there exists
$f=(f_x)$ such that
\begin{eqnarray}
\label{072609}
&&
f_x-\sum_{y}a_{xy}f_y=g_x,\quad x\in\bbZ,\\
&&
\sum_{x}\langle x\rangle^{2k}|f_x|^2<+\infty,\quad \forall \,k=1,2,\ldots.\nonumber
\end{eqnarray}
Observe that the infinite summation range appearing in the first equation is
in fact finite. Indeed, from \eqref{finite-bw} and symmetry it follows that
\begin{equation}
\label{finite-bwa}
a_{x,x-z}=0\mbox{ for all }z>c_n\langle x\rangle^{\ga},
\end{equation}
where
\begin{equation}
\label{c-n}
c_n:=\max\left\{2n,2^{(1+\ga/2)/(1-\ga)}n^{1/(1-\ga)}\right\}.
\end{equation}  
Combining this with \eqref{finite-bw} we conclude that
\begin{equation}
\label{finite-bwb}
a_{x,x+z}=0\mbox{ for all }|z|>c_n\langle x\rangle^{\ga}.
\end{equation}


Furthermore, note that  any $f$ satisfying  conditions
\eqref{072609} belongs to $D(\bar A)$.  Indeed, consider
$f^{(N)}_x:=f_x1_{[|x|\le N]}$. Thanks to the second condition of
\eqref{072609} and \eqref{042609}
we can easily argue that $(I-i A)f^{(N)}\to g$ and $f^{(N)}\to f$, strongly
in $\ell_2(\bbZ)$, as $N\to+\infty$. Since $\bar A$
is the closure of $A$ we conclude that $f\in D(\bar A)$ and 
$(I-i\bar A)f=g$.
Due to the fact that $c_0(\bbZ)$ is dense
in $\ell_2(\bbZ)$ and that the range $R(I-i\bar A)$ is closed we conclude that
$R(I-i\bar A)=\ell_2(\bbZ)$. The proof of the second equality in \eqref{062609}
goes along the same lines. What yet remains to be shown is therefore
\eqref{072609}.

\bigskip

{\em Proof of \eqref{072609}.}  Let 
$$
\chi_N(r):=\left\{
\begin{array}{ll}
r,& |r|\le N,\\
N,&r\ge N,\\
-N,&r\le -N.
\end{array}
\right.
$$  
For a fixed integer  $N$ define $A^{(N)}$ as a bounded, symmetric
operator corresponding to the hermitian matrix whose entries equal
$$
a_{xy}^{(N)}:=\chi_N(a_{xy})1_{[|x-y|\le N]},\quad x,y\in\bbZ.
$$
Given $g\in c_0(\bbZ)$ there is  a (unique) $\tilde f^{(N)}\in\ell_2(\bbZ)$ such that
\begin{equation}
\label{012709}
(I-iA^{(N)}) \tilde f^{(N)}=g.
\end{equation}
We show that for any positive integer $k$ there 
exists a constant $C>0$ such that
\begin{equation}
\label{022709}
\sum_x\langle x\rangle^{2k}|\tilde f^{(N)}_x|^2\le C,\quad N\ge 1.
\end{equation}
Taking this claim for granted (its proof shall be shown momentarily)
we finish the proof of \eqref{072609}.
Using condition \eqref{022709} with any $k>0$ we conclude that the tails of the
infinite sums defining the $\ell_2(\bbZ)$ norms of  $(\tilde f^{(N)})$
are uniformly small in $N$. This proves that the sequence is strongly
precompact in $\ell_2(\bbZ)$, see e.g. Theorem 4.20.1 of \cite{edwards}.
In fact, observe that each $\tilde f^{(N)}\in D(\bar A)$. Indeed, let
$\tilde  f^{(N,M)}:=(\tilde f^{(N)}_x1_{[|x|\le M]})$ for an integer
$M\ge 1$. Of course
$\tilde f^{(N,M)}\in c_0(\bbZ)\subset D(\bar A)$, and it  converges to $\tilde f^{(N)}$ strongly in
$\ell_2(\bbZ)$, as $M\to+\infty$.
On the other hand, from \eqref{022709}  for any $k,\tilde c>0$ 
there exists $C>0$ such that
\begin{equation}
\label{022709b}
\sup_{|y|\ge \tilde c|x|} |\tilde f^{(N)}_y|^2\le \frac{C}{\langle x\rangle^{2k+4}},\quad N\ge 1,\,x\in\bbZ.
\end{equation}
Using \eqref{finite-bwb},  we can estimate 
\begin{eqnarray}
\label{062909}
&&
\sum_x\langle x\rangle^{2k}| (A\tilde  f^{(N,M)})_x|^2
\le \sum_x\langle x\rangle^{2k} \left(\sum_{|y-x|\le c_n\langle
    x\rangle^{\ga}} |a_{xy}| |\tilde f^{(N)}_y|\right)^2\nonumber\\
&&
\le \sum_x\langle x\rangle^{2k}\sup_{|y-x|\le c_n\langle
    x\rangle^{\ga}} |\tilde f^{(N)}_y|^2 \left(\sum_y |a_{xy}|
  \right)^2.
\end{eqnarray}
Since $\ga\in[0,1)$ condition $|y-x|\le c_n\langle
    x\rangle^{\ga}$ implies that there exists $\tilde c>0$ such that
    $|y|\ge \tilde c|x|$ for all $x,y\in\bbZ$.
Thanks to \eqref{022709b} the utmost right hand side of \eqref{062909} can be estimated then by
\begin{equation}
\label{062909a}
 \sum_x\langle x\rangle^{2k}\sup_{|y|\ge \tilde c|x|} |\tilde f^{(N)}_y|^2 \left(\sum_y |a_{xy}|
  \right)^2\le  C\sum_x\langle x\rangle^{-4} \left(\sum_y |a_{xy}|
  \right)^2.
\end{equation}
This together with \eqref{042609} 
imply that there exists $C_1>0$ such that
\begin{equation}
\label{032909}
\sum_x\langle x\rangle^{2k}| (A\tilde  f^{(N,M)})_x|^2
\le C_1\sum_x\langle x\rangle^{-2-2\ga},\quad N,M\ge 1.
\end{equation}
In consequence $(A\tilde  f^{(N,M)})$, $M\ge1$ is strongly precompact
in $\ell_2(\bbZ)$, for a fixed $N$, and
since $\bar A$ is the closure of $A$ we obtain
$\tilde  f^{(N)}\in D(\bar A)$ and
$$
\bar A\tilde  f^{(N)}=\lim_{M\to+\infty}A\tilde  f^{(N,M)}.
$$
 In addition, we also infer that
\begin{equation}
\label{032909a}
(\bar A\tilde  f^{(N)})_x=\sum_ya_{xy}\tilde  f^{(N)}_y,\quad x\in\bbZ
\end{equation}
and that for any $k>0$
 there exists
a constant $C>0$
\begin{equation}
\label{032909n}
\sum_x\langle x\rangle^{2k}| (\bar A\tilde  f^{(N)})_x|^2
\le C,\quad N\ge 1.
\end{equation}
From \eqref{022709} and \eqref{032909n} we conclude that both sequences  $(\tilde  f^{(N)})$ and
$(\bar A\tilde  f^{(N)})$ are strongly precompact in
$\ell_2(\bbZ)$. Choosing a suitable subsequences if necessary we can
assume with no loss of generality that  $\tilde  f^{(N)}\to f$ and
$\bar A\tilde  f^{(N)}$ converges to some $h$ strongly in $\ell_2(\bbZ)$, as $N\to+\infty$. Then $f\in D(\bar A)$
and $\bar Af=h$. In light of \eqref{012709}, to finish the proof of  \eqref{072609}   it suffices to show 
that 
\begin{equation}
\label{042909}
\lim_{N\to+\infty}\|A^{(N)}\tilde  f^{(N)}-\bar A\tilde  f^{(N)}\|_{\ell_2(\bbZ)}=0.
\end{equation}
Estimating as in \eqref{062909} and  \eqref{062909a} we obtain
\begin{eqnarray*}
&&
\sum_x| ((\bar A-A^{(N)})\tilde  f^{(N)})_x|^2
\\
&&
\le \sum_x\sup_{|y|\ge \tilde c|x|} |\tilde f^{(N)}_y|^2 \left(\sum_y
  |a_{xy}-a^{(N)}_{xy}|
  \right)^2 \\
&&
\le C\sum_x\langle x\rangle^{-2k-4} \left(\sum_y
  |a_{xy}-a^{(N)}_{xy}|\right)^2
\to0,
\end{eqnarray*}
as $N\to+\infty$. The passage to the limit  on the utmost right hand side can be argued easily
by virtue of the Lebesgue dominated convergence theorem.
This
ends the proof of \eqref{072609}, modulo the fact that
estimate  \eqref{022709} still requires to be shown. Its proof is an
adaptation to the present case of an
argument used in \cite{sethuraman}, see also Section 2.7.4 of
\cite{klo}. Define a bounded  operator $T:\ell_2(\bbZ)\to\ell_2(\bbZ)$ by the
formula
$Tf_x:=t_xf_x$, where 
\begin{equation}
\label{152709}
t_x:=\langle X\rangle^k1_{[|x|<X]}+\langle x\rangle^k1_{[X\le |x|\le Y]}+\langle Y\rangle^k1_{[Y<|x|]}, \quad x\in\bbZ
\end{equation}
and $0<X<Y$ are some constants  to be determined later on.
Directly from \eqref{152709} it follows that
\begin{equation}
\label{152709a}
|t_x-t_y|\le |\langle x\rangle^k-\langle y\rangle^k|,\quad\forall\,x,y\in\bbZ.
\end{equation}
Applying $T$ to both sides
of \eqref{012709} and taking inner product against $T\tilde f^{(N)}$ on both sides of the aforementioned equation
we conclude that
\begin{eqnarray}
\label{072709}
&&
\|T\tilde f^{(N)}\|_{\ell_2(\bbZ)}^2 +i\langle T\tilde f^{(N)}, [T,A^{(N)}]\tilde f^{(N)}\rangle_{\ell_2(\bbZ)}+i\langle T\tilde f^{(N)}, A^{(N)}T\tilde f^{(N)}\rangle_{\ell_2(\bbZ)}\nonumber\\
&&
=\langle T\tilde f^{(N)}, Tg\rangle_{\ell_2(\bbZ)},
\end{eqnarray}
where $[T,A^{(N)}]:= TA^{(N)}-A^{(N)}T$ is the commutator of $T$ and $A^{(N)}$. 
Thanks to symmetry of  $A^{(N)}$ we have
$$
{\rm Re}\,i\langle T\tilde f^{(N)}, A^{(N)}T\tilde f^{(N)}\rangle_{\ell_2(\bbZ)}=0.
$$
Here ${\rm Re}\,z$ and ${\rm Im}\,z$ denote the real and imaginary parts of a complex number $z$.
Taking the real part of the expressions appearing on both sides of \eqref{072709}  
we obtain 
\begin{eqnarray}
\label{072709a}
&&
\|T\tilde f^{(N)}\|_{\ell_2(\bbZ)}^2 -{\rm Im}\langle T\tilde f^{(N)}, [T,A^{(N)}]\tilde f^{(N)}\rangle_{\ell_2(\bbZ)}\nonumber\\
&&
={\rm Re}\langle T\tilde f^{(N)}, Tg\rangle_{\ell_2(\bbZ)}.
\end{eqnarray}
Note that
$$
[T,A^{(N)}]\tilde f^{(N)}_x=\sum_{y}a^{(N)}_{xy}(t_x-t_y)f_y^{(N)},
$$
therefore 
$$
\langle T\tilde f^{(N)}, [T,A^{(N)}]\tilde f^{(N)}\rangle_{\ell_2(\bbZ)}=\sum_x \sum_{y}a^{(N)}_{yx}(t_x-t_y)t_x\tilde f_x^{(N)}(\tilde f_y^{(N)})^*.
$$
The above expression can be bounded as follows
\begin{equation}
\label{112709}
|\langle T\tilde f^{(N)}, [T,A^{(N)}]\tilde f^{(N)} \rangle_{\ell_2(\bbZ)}|\le \sum_x
\sum_{y}
|a_{xy}| t_x|t_y-t_x||\tilde f_x^{(N)}||\tilde f_y^{(N)}|.
\end{equation}
Applying Young's inequality we can estimate
the right hand side of \eqref{112709} by $I_1+I_2$, where
\begin{eqnarray}
\label{082709}
&&
 I_1:=\frac12\sum_x
\sum_{y}|a_{xy}| t_x|t_y-t_x||\tilde f_x^{(N)}|^2,\\
&&
I_2:=\frac{1}{2}\sum_y
\sum_{x} |a_{xy}| t_x|t_y-t_x||\tilde f_y^{(N)}|^2 \nonumber.
\end{eqnarray} 
Using condition \eqref{finite-bwb}  and the fact that  $t_x$ is constant for $|x|\le X$, or $Y\le |x|$ we conclude that
$$
 I_1\le \frac12\sum_{x=\bar X}^{\bar Y}
\sum_{y=x-c_n\langle x\rangle^{\ga}}^{x+c_n\langle x\rangle^{\ga}}
|a_{xy}| t_x|t_x-t_y| |\tilde f_x^{(N)}|^2,
$$
where
$\bar X:=X-c_n\langle X\rangle^{\ga}$ and $\bar Y:=Y+c_n\langle
Y\rangle^{\ga}$. Thanks to \eqref{152709a} we can estimate
\begin{equation}
\label{103009}
 I_1\le \frac12\sum_{x=\bar X}^{\bar Y}  t_x^2  |\tilde f_x^{(N)}|^2 \left\{\frac{\langle x\rangle^k}{t_x}
\sum_{y=x-c_n\langle x\rangle^{\ga}}^{x+c_n\langle x\rangle^{\ga}}
|a_{xy}|\left|\frac{\langle y\rangle^k}{\langle x\rangle^k}-1\right|\right\}.
\end{equation}
 Observe that 
 \begin{equation}
 \label{010610}
 \frac{\langle x\rangle^k}{t_x}\le (c_n+1)^{k/2},\quad  \mbox{for }\bar X\le x\le \bar Y.
 \end{equation}
Choose an arbitrary $\delta\in(0,1)$. Note that for $|m|\le n\langle
x\rangle^{\ga}$ there exist constants $C,C'>0$ such that
\begin{equation}
\label{102709}
\left|\left(\frac{1+
    (x+m)^2}{1+
    x^2}\right)^{k/2}-1\right|\le C\frac{|xm|+m^2/2}{1+x^2}\le
C'(\langle x\rangle^{\ga-1}+\langle x\rangle^{2\ga-2})
\end{equation}
for all $x$.
Combining \eqref{102709} with \eqref{042609} we conclude that for any $\delta\in(0,1)$ there exists $X$, depending on parameters
$n,\delta,k,\ga$, such that 
\begin{equation}
\label{142709}
\sum_{y=x-c_n\langle x\rangle^{\ga}}^{x+c_n\langle x\rangle^{\ga}} |a_{xy}|\left|\frac{\langle
    y\rangle^k}{\langle x\rangle^k}-1\right|<\frac{\delta}{(c_n+1)^{k/2}},\quad \mbox{ for
}|x|\ge \bar X.
\end{equation}
This together with \eqref{010610} imply that
\begin{equation}
\label{172709}
I_1\le \frac{\delta}{2} \|T\tilde f^{(N)}\|_{\ell_2(\bbZ)}^2.
\end{equation}
Likewise
$$
 I_2\le \frac12\sum_{y=\bar X}^{\bar Y}  t_y^2  |\tilde f_y^{(N)}|^2 \sum_{x=y-c_n\langle y\rangle^{\ga}}^{y+c_n\langle y\rangle^{\ga}}
\frac{ t_x}{t_y^2}
|a_{xy}|\left|\langle x\rangle^k-\langle y\rangle^k\right|.
$$
Since
$$
\frac{1}{C_*}\le \frac{t_x}{\langle x\rangle^k}\le C_*,\quad x=\bar
X- c_n\langle \bar X\rangle^{\ga},\ldots,\bar
Y+ c_n\langle \bar Y\rangle^{\ga}
$$
for some constant $C_*>0$ that depends only on $n,\ga,k$, we conclude that
$$
 I_2\le \frac{C_*^3}{2}\sum_{y=\bar X}^{\bar Y}  t_y^2  |\tilde  f_y^{(N)}|^2 \sum_{x=y-c_n\langle y\rangle^{\ga}}^{y+c_n\langle y\rangle^{\ga}}|a_{xy}|
\frac{\langle
    x\rangle^k}{\langle y\rangle^k}\left|\frac{\langle
    x\rangle^k}{\langle y\rangle^k}-1\right|.
$$
Choose $X$ sufficiently large so that
$$
\sum_{x=y-c_n\langle y\rangle^{\ga}}^{y+c_n\langle y\rangle^{\ga}} |a_{xy}| \frac{\langle
    x\rangle^k}{\langle y\rangle^k}\left|\frac{\langle
    x\rangle^k}{\langle y\rangle^k}-1\right|<\frac{\delta}{C_*^3},\quad \mbox{ for }
|y|\ge \bar X.
$$
As a result we conclude
\begin{equation}
\label{172709a}
I_2\le \frac{\delta}{2} \|T\tilde f^{(N)}\|_{\ell_2(\bbZ)}^2.
\end{equation}
Combining this with \eqref{172709} 
we have
$$
|\langle T\tilde f^{(N)}, [T,A^{(N)}]\tilde f^{(N)} \rangle_{\ell_2(\bbZ)}|\le \delta\sum_xt_x^2| \tilde f_x^{(N)}|^2=\delta\|T \tilde f^{(N)}\|_{\ell_2}^2.
$$
Going back to \eqref{072709a} we obtain
\begin{eqnarray*}
&&
\|T\tilde f^{(N)}\|_{\ell_2(\bbZ)}^2 \le  |\langle T\tilde f^{(N)},
Tg\rangle_{\ell_2(\bbZ)}|+ |\langle T\tilde f^{(N)}, [T,A^{(N)}]\tilde f^{(N)} \rangle_{\ell_2(\bbZ)}|\\
&&
\le |\langle T\tilde f^{(N)},
Tg\rangle_{\ell_2(\bbZ)}|+\delta \|T\tilde f^{(N)}\|_{\ell_2(\bbZ)}^2,
\end{eqnarray*}
therefore
\begin{equation}
\label{023009}
\|T\tilde f^{(N)}\|_{\ell_2(\bbZ)}\le \frac{\|Tg\|_{\ell_2(\bbZ)}}{1-\delta}.
\end{equation}
Now, let parameter $Y$, appearing in the definition of the operator $T$, tend to
infinity. Since $\|Tg\|_{\ell_2(\bbZ)}$ remains constant, starting with some sufficiently large $Y$ (as $g\in c_0(\bbZ)$) we infer that \eqref{023009} implies \eqref{022709}. This ends the proof
of \eqref{072609} finishing also the proof of Theorem \ref{main}.

\bigskip

\section{Proof of Theorem \ref{cor-main1}}

The argument used in the proof of  Theorem \ref{main} can be
adapted  to the present case, provided we are able to show that \eqref{072609}
holds for some fixed $k_0>2$. To prove this fact we repeat with no changes,
except replacing $c_n$ by $n$ (maintaining the notation from the previous section), the calculations made 
between \eqref{152709} and \eqref{103009}. 
Instead of \eqref{102709} we write that for some constant $C>0$
\begin{equation}
\label{102709a}
\left|\left(\frac{1+
    (x+m)^2}{1+
    x^2}\right)^{k_0/2}-1\right|\le \frac{C}{\langle x\rangle},\quad
\forall\,x\in\bbZ, \,|m|\le n.
\end{equation}
Using the above together with  \eqref{042609} we conclude that for any
$\eps>0$ there exists $X$ (appearing in the definition of operator $T$), depending on 
$n$, such that 
\begin{equation}
\label{142709a}
\sum_{y=x-n}^{x+n} |a_{xy}|\left|\frac{\langle
    y\rangle^{k_0}}{\langle x\rangle^{k_0}}-1\right|\le c_*+\eps,\quad \mbox{ for
}|x|\ge \bar X
\end{equation}
and
$$
\sum_{x=y-n}^{y+n} |a_{xy}| \frac{\langle
    x\rangle^{k_0}}{\langle y\rangle^{k_0}}\left|\frac{\langle
    x\rangle^{k_0}}{\langle y\rangle^{k_0}}-1\right|\le c_*+\eps,\quad \mbox{ for }
|y|\ge \bar X.
$$
Here $\bar X:=X-n$.
This leads to an estimate 
\begin{equation}
\label{172709s}
I_1+I_2\le \frac12 (c_*+\eps)\left[(n+1)^{k_0/2} +C_*^3\right]\|T\tilde f^{(N)}\|_{\ell_2(\bbZ)}^2.
\end{equation}
Choosing $c_*$ and $\eps>0$ in such a way that 
$$
\delta:=\frac12 (c_*+\eps)\left[(n+1)^{k_0/2} +C_*^3\right]<1
$$
we can still claim \eqref{023009}. This allows us to conclude  \eqref{072609}, which ends the proof of the theorem.

\end{document}